 \UseRawInputEncoding
\documentclass[sn-mathphys,Numbered]{sn-jnl}

\usepackage{graphicx}%
\usepackage{multirow}%
\usepackage{amsmath,amssymb,amsfonts}%
\usepackage{amsthm}%
\usepackage{mathrsfs}%
\usepackage[title]{appendix}%
\usepackage{xcolor}%
\usepackage{textcomp}%
\usepackage{manyfoot}%
\usepackage{booktabs}%
\usepackage{algorithm}%
\usepackage{algorithmicx}%
\usepackage{algpseudocode}%
\usepackage{listings}%
\usepackage{latexsym, amsmath, xcolor, graphicx, amssymb, color,booktabs,mathrsfs,listings,subcaption,float,stackrel,amsthm,mathdots}
\usepackage{multirow}
\usepackage{xcolor, soul}
\sethlcolor{yellow}
\usepackage{mathtools}
\usepackage{bm}
\usepackage{lineno}

\newcommand{\argmin}{\mathop{\mathrm{arg\,min}}}
\begin{document}

\title[Estimation of Parameter Distributions for Reaction-Diffusion Equations with Competition using Aggregate Spatiotemporal Data]{Estimation of Parameter Distributions for Reaction-Diffusion Equations with Competition using Aggregate Spatiotemporal Data}


\author[1,2]{\fnm{Kyle} \sur{Nguyen}}\email{kcnguye2@ncsu.edu}

\author[3]{\fnm{Erica M.} \sur{Rutter}}\email{erutter2@ucmerced.edu}

\author*[2,4]{\fnm{Kevin} \sur{Flores}}\email{kbflores@ncsu.edu}

\affil[1]{\orgdiv{Biomathematics Graduate Program}, \orgname{North Carolina State University}, \orgaddress{\city{Raleigh}, \state{NC}, \country{USA}}}

\affil[2]{\orgdiv{Center for Research in Scientific Computation}, \orgname{North Carolina State University}, \orgaddress{\city{Raleigh}, \state{NC}, \country{USA}}}

\affil[3]{\orgdiv{Department of Applied Mathematics}, \orgname{University of California, Merced}, \orgaddress{\city{Merced}, \state{CA}, \country{USA}}}

\affil[4]{\orgdiv{Department of Mathematics}, \orgname{North Carolina State University}, \orgaddress{\city{Raleigh}, \state{NC}, \country{USA}}}


\abstract{Reaction diffusion equations have been used to model a wide range of biological phenomenon related to population spread and proliferation from ecology to cancer. It is commonly assumed that individuals in a population have homogeneous diffusion and growth rates, however, this assumption can be inaccurate when the population is intrinsically divided into many distinct subpopulations that compete with each other. In previous work, the task of inferring the degree of phenotypic heterogeneity between subpopulations from total population density has been performed within a framework that combines parameter distribution estimation with reaction-diffusion models. Here, we extend this approach so that it is compatible with reaction-diffusion models that include competition between subpopulations. We use a reaction-diffusion model of Glioblastoma multiforme, an aggressive type of brain cancer, to test our approach on simulated data that are similar to measurements that could be collected in practice. We use Prokhorov metric framework and convert the reaction-diffusion model to a random differential equation model to estimate joint distributions of diffusion and growth rates among heterogeneous subpopulations. We then compare the new random differential equation model performance against other partial differential equation models' performance. We find that the random differential equation is more capable at predicting the cell density compared to other models while being more time efficient. Finally, we use $k$-means clustering to predict the number of subpopulations based on the recovered distributions.
}




\maketitle
\section{Introduction}
\label{introduction}
The importance of including phenotypic heterogeneity in partial differential equation models (PDE) of spatially diffusing or phenotypically structured populations, including reaction-diffusion, advection-diffusion, and Sinko-Streifer type models, has been emphasized previously by Banks and Knunisch \cite{banks_estimation_1989} and discussed further in \cite{banks_quantifying_2007}. The traditional approach to incorporating phenotypic heterogeneity is to parameterize phenotypic differences by using several PDEs, i.e., one per phenotype \cite{matsiaka2019continuum,hatzikirou_go_2012}. For example, a reaction-diffusion PDE used to model a population that consists of two subpopulations, one that predominately diffuses and another that proliferates, can be extended to a model with two reaction-diffusion PDEs in which the first PDE has a high diffusion and low proliferation rate while the second PDE has a low diffusion and high proliferation rate. An alternative to this traditional approach to incorporating phenotypic heterogeneity into PDE models is to consider parameter coefficients as functional coefficients. This alternative approach has been used in previous modeling studies \cite{banks_parameter_1983, banks_modeling_1985} in which parameters are assumed to be spatially or temporally dependent. However, in many biological settings, data are presented as aggregate data, in which the observed variable(s) is the total population level data \cite{banks_estimation_1999, banks_estimating_2018,rutter_estimating_2018,schacht_estimation_2019,banks_parameter_2020}, i.e., the density of each subpopulation is not individually measured. The inverse problem to recover the distribution of each subpopulation from aggregate data is an ongoing area of research \cite{meyers2021koopman,hatzikirou2021novel}. An inverse problem methodology called the Prokhorov metric framework (PrMF) was previously developed to infer subpopulation heterogeneity from aggregate data by considering parameter coefficients in PDE models as distributed probability functions \cite{banks_functional_2012, banks_modeling_2014}. The PrMF has previously been applied to aggregate temporal data \cite{banks_comparison_2007,banks_estimating_2018,sirlanci_deconvolving_2018,sirlanci_applying_2019,sirlanci2019estimating, schacht_estimation_2019,banks_parameter_2020} and spatiotemporal data \cite{banks_estimation_1999,rutter_estimating_2018}. 

A recent review of the PrMF methodology and its application to mathematical models fit to aggregate data can be found in \cite{banks_parameter_2020}. Performing an inverse problem with the PrMF relies on interpreting parameters in PDE models as distributions of a random differential equation (RanDE), in contrast to the traditional approach of assuming them to be point estimates of a deterministic differential equation. The RanDE is fit to aggregate data by taking an expectation of a forward solution of the RanDE over an unknown parameter distribution. Thereby, the inverse problem of inferring subpopulation heterogeneity is solved by identifying the unknown parameter distribution; i.e., choosing a single distribution among a given family of distributions that best fits the observed data. For example, for the RanDE version of the classic Sinko-Streifer model, advection and mortality distributions can be estimated \cite{banks_parameter_1983,banks_modeling_1985,banks_estimation_1999,banks_comparison_2007,banks_functional_2012,banks_modeling_2014,banks_estimating_2018,schacht_estimation_2019}. Similarly, the Fisher Kolmogorov–Petrovsky–Piskunov (Fisher-KPP) equation \cite{fisher_wave_1937} can be interpreted as a RanDE when the diffusion or growth rates are assumed to be parameter distributions \cite{rutter_estimating_2018,banks_parameter_2020}. Computational methods that frame parameter distribution estimation as a least squares optimization problem as well as the corresponding theoretical underpinnings for convergence are reviewed in \cite{banks_parameter_2020}; additional details for implementation of the PrMF are exemplified with the Fisher-KPP equation in Section 2.

In this work, we focus on extending the PrMF for application to the Fisher-KPP equation \cite{fisher_wave_1937}, a reaction-diffusion equation that can be used to describe the spatial and temporal growth and spread of a population. This equation has been found to be relevant in a wide range of biological settings \cite{murray_mathematical_2002}. The Fisher-KPP equation has been used to model the spread of species in ecology \cite{skellam_random_1951,van_den_bosch_analysing_1992,shigesada_modeling_1995,steele_modelling_1998,reise_item_2009,kuehn_warning_2013}; in epidemiology \cite{mollison_dependence_1991,hethcote_mathematics_2000}; cell migration in wound healing \cite{cai_multi-scale_2007,tremel_cell_2009,habbal_assessing_2014,nardini_modeling_2016,nardini_investigation_2018}; and tumor growth \cite{swanson_quantitative_2000,murray_mathematical_2002,wang_prognostic_2009,neal_response_2013,baldock_patient-specific_2014,le_mri_2016,rutter_mathematical_2017,hawkins-daarud_quantifying_2019}. The specific biological application of the Fisher-KPP equation we consider in this paper is towards modeling the spatial spread and growth of an aggressive type of brain cancer called Glioblastoma multiforme (GBM). The Fisher-KPP equation has been extensively utilized in the literature for modeling GBM \cite{swanson_quantitative_2000,murray_mathematical_2002,wang_prognostic_2009,neal_response_2013,baldock_patient-specific_2014,le_mri_2016,rutter_mathematical_2017,hawkins-daarud_quantifying_2019} and previous work has shown that model predictions are clinically relevant \cite{wang_prognostic_2009,neal_response_2013,baldock_patient-specific_2014}. The single PDE Fisher-KPP model assumes a homogeneous growth and motility phenotype among all tumor cells, however, Stein et al. \cite{stein_mathematical_2007} suggested that cells at the core and cells at the rim exhibit different behaviors. In particular, the ``go or grow" hypothesis \cite{hatzikirou_go_2012} suggests that there are two types of GBM cells: proliferative cells and invasive cells, in which proliferative cells grow faster and are less motile than invasive cells \cite{hatzikirou_go_2012}. Therefore, this introduces the need to model proliferative and invasive cells separately. Since the introduction of the ``go or grow" hypothesis, many mathematical models have been created to model this behavior, normally as a system of partial differential equations in which the parameters governing each population are different  \cite{martinez-gonzalez_hypoxic_2012,pham_density-dependent_2012,stepien_traveling_2018,zhigun_strongly_2018,tursynkozha_traveling_2023}. 

To expand further from the ``go or grow'' hypothesis on spatio-temporal data from models with only 2 subpopulations to a model that uses parameter distributions to define subpopulations, Rutter et al. \cite{rutter_estimating_2018} utilized the PrMF to perform an inverse problem for estimating growth and diffusion rate distributions for a Fisher-KPP model of GBM. Unlike the previous coupled PDE models which assumed two distinct subpopulations, the RanDE is flexible enough to model any number of subpopulations. The proposed method was used to estimate the diffusion coefficient $D$, and the proliferation rate $\rho$ as independent distributions using the PrMF \cite{rutter_estimating_2018} on generated synthetic data with heterogeneous subpopulations. The metric was evaluated on a variety of probability distributions including independent bigaussian distributions, i.e. distributions for the ``grow or go" hypothesis, for $D$ and $\rho$ \cite{rutter_estimating_2018}. While the authors were able to recover the underlying distributions \cite{rutter_estimating_2018}, there are some limitations. One of the limitations is that the distributions of $D$ and $\rho$ are assumed to be independent rather than a joint distribution. In addition, the competition between subpopulations was neglected within the proposed model.

In this work, we propose an improved inverse problem modeling method applying the PrMF and built upon the work of Rutter et al. \cite{rutter_estimating_2018}. We note that, while the application of the PrMF is specific to the Fisher-KPP equation in this paper, one of the primary contributions of our work is to exemplify how the PrMF could be applied to the general scenario in which a differential equation model with competition among heterogeneous subpopulations is being fit to aggregate data. Another key contribution is the first implementation of the PrMF to a joint distribution over two parameters that are treated as random variables. Previous efforts utilizing the PrMF have only estimated a single parameter distribution or treated a joint distribution as two independent distributions. Estimating a joint distribution is more accurate in many biological scenarios and we exemplify it here for the case of the Fisher-KPP model for GBM. In particular, modeling subpopulations to investigate the go or grow hypothesis would require a subpopulation that \textit{grows} and a subpopululation \textit{goes}. Thus, one subpopulation should have a distribution that is centered around a high diffusion rate and low proliferation rate while the other subpopulation should have low diffusion and high proliferation. The previous implementation of the PrMF carried out in \cite{rutter_estimating_2018} does not allow for such a scenario. Similar to the work of Rutter et al. \cite{rutter_estimating_2018}, we evaluate our framework on noisy synthetic data with heterogeneous subpopulations. Particularly, we generate the data from a RanDE model with a mixture of gaussian distributions for $D$ and $\rho$. This means $D$ and $\rho$ are assumed to be random variables. In this work, their distributions are approximated by discrete nodes. Each pair of $D$ and $\rho$ has an associated weight describing a discrete probability density function (see section. \ref{discrete_approx_PMF} for more details). We then compare the RanDE model performance against the classical reaction-diffusion models with 2, 4, and 6 subpopulations, i.e. models with 2, 4, and 6 reaction-diffusion PDEs that includes competition. However, we rely on the traditional inverse problem approach to perform point-wise parameter estimation on these classical PDEs. On the other hand, we apply the PrMF on the RanDE model to perform distribution estimation on the RanDE model. The models are then compared in terms of fitting, predicting, and computation time. Finally, we use $k$-means clustering, a simple unsupervised machine learning method, to recover the number of subpopulations and their cluster centers. Finally, we briefly describe the use of $k$-means clustering for recovering the number of subpopulations in section \ref{kmean}.

\section{Methods}\label{methods}
In this section we detail all methodologies used to compare inference of subpopulation heterogeneity using PrMF on RanDE model and the traditional method using a set of coupled PDEs. In section \ref{ranDE_model}, we outline the mathematical models. This is followed by the introduction of the PrMF approach for joint distribution estimation in section \ref{PMF}. In section \ref{appendix_A}, we briefly describe the traditional inverse problem approach on a set of coupled PDEs to perform point-wise parameter estimation. We follow with our data generation process in section \ref{data}.

\subsection{Mathematical models} \label{ranDE_model}
The Fisher-KPP equation is a reaction-diffusion PDE with a single population:

\begin{equation}
    \frac{\partial u(x,t)}{\partial t} = D \frac{\partial^2 u(x,t)}{\partial x^2} + \rho u(x,t) \left[1 - \frac{u(x,t)}{K} \right]
    \label{fisherkpp}
\end{equation}
where $u(x,t)$ is the aggregated cell density at the spatial coordinate $x$ and temporal coordinate $t$. In addition, $D$ is the diffusion coefficient, corresponding to the invasiveness of the cells to the surrounding areas. $\rho$ is the growth rate and $K$ is the carrying capacity for the cell density. For a normalized Fisher-KPP equation, we set $K=1$ and obtain the following equation:

\begin{equation}
    \frac{\partial u(x,t)}{\partial t} = D \frac{\partial^2 u(x,t)}{\partial x^2} + \rho u(x,t) \left[1 - u(x,t) \right]
    \label{normalized_fisherkpp}
\end{equation}

The normalized Fisher-KPP equation is often used as the simplest model to describe the spatial and temporal growth and spread of GBM \cite{swanson_quantitative_2000,murray_mathematical_2002,wang_prognostic_2009,neal_response_2013,baldock_patient-specific_2014,le_mri_2016,rutter_mathematical_2017,hawkins-daarud_quantifying_2019}. However, this model with a single population often fails to describe the heterogeneous behavior in GBM tumor, where the invasive cells near the rim exhibit different behavior compared to the proliferative cells at the center \cite{stein_mathematical_2007}. In particular, the invasive cells have the tendency to migrate to the surrounding areas, while the proliferative cells often grow faster compared to the invasive cells \cite{stein_mathematical_2007}. This leads to the suggestion to model invasive cells and proliferative cells separately. The ``go or grow" hypothesis is then introduced with switching functions from proliferative cells to invasive cells and vice versa \cite{hatzikirou_go_2012}. Later, Stepien et al. introduced a general form of the ``go or grow" model for spatial-temporal data with switching functions \cite{stepien_traveling_2018}.

Expanding from the ``go or grow" hypothesis, Rutter et al. introduced a RanDE version of the Fisher-KPP equation \cite{rutter_estimating_2018}. In a RanDE model, the parameters take the form of distributions rather than point estimates in deterministic differential equations. The RanDE for the Fisher-KPP equation introduced by Rutter et al.  is described as follows \cite{rutter_estimating_2018}:

\begin{equation}
    \frac{\partial c(x,t;\boldsymbol{D},\boldsymbol{\rho})}{\partial t} = \boldsymbol{D} \frac{\partial^2 c (x,t;\boldsymbol{D},\boldsymbol{\rho})}{\partial x^2} + \boldsymbol{\rho} c(x,t;\boldsymbol{D}, \boldsymbol{\rho})\left[1 - c(x,t;\boldsymbol{D},\boldsymbol{\rho})\right]
    \label{rutter_subpop}
\end{equation}
where $\boldsymbol{D}$ and $\boldsymbol{\rho}$ are the random variables for the diffusion and growth rates on a compact space, $\Omega=\Omega_{D} \times \Omega_\rho$. $c(x,t;\boldsymbol{D},\boldsymbol{\rho})$ is the spatiotemporal cell density for the phenotype corresponding to random variables $\boldsymbol{D}$ and $\boldsymbol{\rho}$. Assuming the probability measure for $\boldsymbol{D}$ and $\boldsymbol{\rho}$ is $P(\boldsymbol{D},\boldsymbol{\rho})$, then the aggregated population, $u(x,t)$ is computed as follows \cite{rutter_estimating_2018}:

\begin{equation}
    u(x,t;P) = \int_{\Omega} c(x,t;\boldsymbol{D},\boldsymbol{\rho}) dP(\boldsymbol{D},\boldsymbol{\rho})
    \label{aggregated_u}
\end{equation}

One limitation of Eq. \ref{rutter_subpop} is that intratumor competition is neglected. Therefore, to incorporate the intratumor competition between phenotype, we propose the following model:

\begin{equation}
\begin{split}
    \frac{\partial c(x,t;\boldsymbol{D},\boldsymbol{\rho})}{\partial t} &= \boldsymbol{D} \frac{\partial^2 c (x,t;\boldsymbol{D},\boldsymbol{\rho})}{\partial x^2}\\
    & + \boldsymbol{\rho} c(x,t;\boldsymbol{D}, \boldsymbol{\rho})\left[1 - \int_{\Omega} \alpha(\boldsymbol{D},\boldsymbol{\rho}) c(x,t;\boldsymbol{D},\boldsymbol{\rho}) dP(\boldsymbol{D},\boldsymbol{\rho}) \right]
    \label{new_subpop}
\end{split}
\end{equation}
with $\alpha(\boldsymbol{D},\boldsymbol{\rho})$ representing the competitive advantage for the phenotype corresponding to the random variables $\boldsymbol{D}$ and $\boldsymbol{\rho}$. For simplification, we further assume that all phenotype have the same competitive advantage and let $\alpha(\boldsymbol{D},\boldsymbol{\rho}) = 1$ for all $(\boldsymbol{D},\boldsymbol{\rho})$. By substituting Eq. \ref{aggregated_u} into Eq. \ref{new_subpop}, we have:

\begin{equation}
    \frac{\partial c(x,t;\boldsymbol{D},\boldsymbol{\rho})}{\partial t} = \boldsymbol{D} \frac{\partial^2 c (x,t;\boldsymbol{D},\boldsymbol{\rho})}{\partial x^2} + \boldsymbol{\rho} c(x,t;\boldsymbol{D}, \boldsymbol{\rho})\left[1 - u(x,t;P) \right]
    \label{new_subpop2}
\end{equation}
where $u(t,x;P)$ is provided in Eq. \ref{aggregated_u}.

To summarize the model, the aggregated population contains multiple phenotypes. Thus, the aggregate data defined in Equation (4) as the sum of all subpopulations among all phenotypes described by the random variables $\boldsymbol{D}$ and $\boldsymbol{\rho}$ is what is used for data fitting. In this study, we investigate data generated from a joint distribution over $\boldsymbol{D}$ and $\boldsymbol{\rho}$ that describes two subpopulations consisting of proliferative (grow, high $\boldsymbol{\rho}$ and low $\boldsymbol{D}$) and invasive (go, high $\boldsymbol{D}$ and low $\boldsymbol{\rho}$) phenotypes. Importantly, our method does not rely on observations of each respective subpopulation, and instead is able to use the type of information that can be measured in practice, i.e., the total population density. In this study, we compare the performance in term of fitting and predicting between the RanDE model (Eqs. \ref{aggregated_u} and \ref{new_subpop2}) and classical reaction-diffusion competition models that have 2, 4, or 6 reaction-diffusion PDEs that describe populations with 2, 4, and 6 subpopulations respectively. Details for these models can be found in section \ref{appendix_A}. For simplicity, we denote these classical models as 2-, 4-, and 6-PDEs models. Note that Eqs. \ref{aggregated_u} and \ref{new_subpop2} can be viewed as a generalized competition model.

\subsection{Prokhorov metric framework} \label{PMF}

To perform an inverse problem on a random differential equation model, we rely on the Prokhorov Metric Framework (PrMF). The method was developed in \cite{banks_functional_2012} and completed in \cite{banks_modeling_2014}. The framework has been applied to various biological problems with aggregated data. For example, PrMF was applied to estimate the distribution of growth rate for mosquitofish \cite{banks_estimation_1999,banks_comparison_2007} and for prions \cite{banks_estimating_2018}. An inverse problem using PrMF were also used to estimate the distribution rate of T-cells leaving blood to the tumor in chimeric antigen receptor therapies model \cite{schacht_estimation_2019,banks_parameter_2020}. The metric also has been applied to non-biological problems \cite{white_using_2022,white_modeling_2020}. To expand from the ``go or grow" hypothesis, Rutter et al. proposed to perform an inverse problem to estimate the distributions diffusion $\boldsymbol{D}$ and growth rate $\boldsymbol{\rho}$ for the random differential equation version of the reaction-diffusion equation (Fisher-KPP equation) \cite{rutter_estimating_2018}. However, one of the limitations of this previous study is that the distributions for $\boldsymbol{D}$ and $\boldsymbol{\rho}$ are assumed to be independent.

In this study, we generalize the inverse problem approach using PrMF by assuming the interested distribution to be joint distributions of $\boldsymbol{D}$ and $\boldsymbol{\rho}$ on aggregated cell density data. Particularly, we approximate the probability measure $P(\boldsymbol{D},\boldsymbol{\rho})$ for each pair of $(\boldsymbol{D},\boldsymbol{\rho})$. In other words, we estimate the contribution density of the subpopulation corresponding to $(\boldsymbol{D},\boldsymbol{\rho})$ towards the aggregated population data. There are two different methods to approximate $P(\boldsymbol{D},\boldsymbol{\rho})$: delta functions method and spline functions methods \cite{banks_comparison_2007,rutter_estimating_2018}. The method based in delta functions recovers a discrete estimation of the distribution while the method using splines recovers a continuous estimation of the distribution. In this study, we only use the delta functions method. For details about the spline functions, we refer the readers to previous studies \cite{banks_comparison_2007,rutter_estimating_2018}.

\subsubsection{Discrete approximation} \label{discrete_approx_PMF}

Let $M_{D}$ and $M_{\rho}$ be the numbers of evenly spaced sampled nodes for $D$ and $\rho$. We denote the discretized version of $(\boldsymbol{D},\boldsymbol{\rho})$ as $(D_i,\rho_i)$, where $D_i \in \left\{ \min_D=D_0, \dots, D_{M_D}=\max_D \right\}$ and $\rho_i \in \left\{ \min_\rho=\rho_0, \dots, \rho_{M_\rho}=\max_\rho \right\}$. We further denote $P(\boldsymbol{D},\boldsymbol{\rho})$ as $\boldsymbol{w}$. Additionally, $M=M_D \times M_\rho$ is the number of cells being used for discretization of the joint probability density function described by $\boldsymbol{w} = \{w_i\}_{i=1}^M$. $M$ also represents the number of phenotypes in the aggregated population. Using discrete approximation, we simplify Eqs. \ref{aggregated_u} and \ref{new_subpop2} to:

\begin{equation}
    u \left( x,t;\boldsymbol{w} \right) = \sum_{i}^{M} w_i c(x,t;D_i,\rho_i)
    \label{aggregated_discrete_u}
\end{equation}
and

\begin{equation}
    \frac{\partial c(x,t;D_i,\rho_i)}{\partial t} = D_i \frac{\partial^2  c(x,t;D_i,
    \rho_i)}{\partial x^2} + \rho_i c(x,t;D_i,\rho_i)\left[1 - u(x,t;\boldsymbol{w}) \right]
    \label{new_discrete_subpop}
\end{equation}
where $w_i \geq 0$ is the discrete weight associated with the $i$th phenotype. In other words, $w_i$ is the probability measure for $(D_i,\rho_i)$. This leads to the constraint $\sum_{i}^{M}w_i = 1$. Finally, we define the finite dimensional approximation to $P(\Omega)$ with $M$ cells as:

\begin{equation}
    P^{M}(\Omega) = \Bigl\{ \boldsymbol{w} \in P( \Omega ) \vert \boldsymbol{w} = \sum_{i}^M w_i \delta_{D_i,\rho_i}, \; \sum_{i}^M w_i = 1 \text{ and } w_i \geq 0\Bigl\}
\end{equation}
where $\delta_{D_i,\rho_i}$ is the delta function with an atom at $(D_i,\rho_i)$. 

\subsubsection{Forward solve phenotype equations with competition} \label{forwardsolve_PMF}

To perform parameter estimation using the PrMF method, we need to approximate the phenotype cell densities, $c(x,t;D_i,\rho_i)$, in Eq. \ref{new_discrete_subpop}. Therefore, we forward solve for the solutions of Eq. \ref{new_discrete_subpop} for different $(D_i,\rho_i)$ with $i=1, \dots,M$. Unlike the RanDE model without competition developed by Rutter et al. \cite{rutter_estimating_2018}, the new model (Eqs. \ref{aggregated_discrete_u} and \ref{new_discrete_subpop}) requires prior information about the aggregated population, $u (x,t;\left\{ w\right\}_{i=1}^M )$. However, we can assume that the aggregated population is the observed data, $u^o(t,x)$. Therefore, we simplify Eq. \ref{new_discrete_subpop} to:

\begin{equation}
    \frac{\partial c(x,t;D_i,\rho_i)}{\partial t} = D_i \frac{\partial^2 c(x,t;D_i,
    \rho_i)}{\partial x^2} + \rho_i c(x,t;D_i,\rho_i)\left[1 - u^o(x,t) \right]
    \label{new_discrete_subpop_data}
\end{equation}

For each phenotype, we use central differencing to approximate the spatial derivatives in Eq. \ref{new_discrete_subpop_data}. Next, we numerically solve the model using the Explicit Runge-Kutta (RK45) method as the integration method, implemented  in \verb|scipy| as \verb|integrate| . 


\subsubsection{Estimation of discrete weights} \label{estimate_PMF_weights}

After we forward solve for the solution of phenotype, we optimize the parameters, i.e. the discrete weights $\boldsymbol{w}$, by minimizing the sum of squares ($\textit{SSE}$). We estimate the parameter vector $\boldsymbol{\hat{w}}$ in the RanDE model such that:

\begin{equation}
    \boldsymbol{\hat{w}} = \argmin_{\boldsymbol{w} \in P^M(\Omega)} \sum_{j,k}^{N_t, N_x} \left[ u_{j,k}^o - \hat{u}(t_j,x_k;\boldsymbol{w}) \right]^2
    \label{argmin_cont}
\end{equation}
where $u_{j,k}^o$ and $\hat{u}(t_j,x_k;\boldsymbol{w})$, respectively, are the observed and simulated values for the cell density at the $j$th time point and $k$th spatial point. $N_t$ is the number of time points while $N_x$ is the number of spatial points. By substituting Eq. \ref{aggregated_discrete_u} into Eq. \ref{argmin_cont}, the equation becomes:

\begin{equation}
    \boldsymbol{\hat{w}} = \argmin_{\boldsymbol{w} \in P^M(\Omega)} \sum_{j,k}^{N_t, N_x} \left[ u_{j,k}^o - \sum_{i}^{M} w_i c(x,t;D_i,\rho_i) \right]^2
    \label{argmin_discrete}
\end{equation}

We use \verb|scipy| package with the built-in function \verb|optimize.minimize| to minimize the $\textit{SSE}$ with bounds and constraints over a fitting time interval. By setting the default option, we let the function choose the Sequential Least SQuares Programming (SLSQP) as the bound and constrained optimization method \cite{kraft_software_1988}.

\subsection{Traditional inverse problem approach for coupled PDEs models} \label{appendix_A}
The competition model with $M$ partial differential equations is described as follows: 

\begin{equation}
    \frac{\partial c(x,t;D_i,\rho_i)}{\partial t} = D_i \frac{\partial^2  c(x,t;D_i,
    \rho_i)}{\partial x^2} + \rho_i c(x,t;D_i,\rho_i)\left[1 - u(x,t) \right]
\end{equation}
with

\begin{equation}
    u \left( x,t \right) = \sum_{i}^{M} w_i c(x,t;D_i,\rho_i),
\end{equation}
where $u(x,t)$ is the aggregated population. $D_i$ and $\rho_i$ are the diffusion coefficient and growth rate of the $i$th subpopulation, respectively. The weight for each subpopulation is denoted as $w_i$. In this work, for the 2-PDE, 4-PDE, and 6-PDE models, $M=2,4,$ and 6, respectively.
During optimization step, we estimate the point-wise parameter vector, $\bold{\hat{q}} = \left(D_1, \dots, D_M, \rho_1, \dots, \rho_M, w_1, \dots, w_M \right)$, by minimizing the $\textit{SSE}$ as follows:

\begin{equation}
    \boldsymbol{\hat{q}} = \argmin_{\boldsymbol{q} \in \bold{Q}} \sum_{j,k}^{N_t, N_x} \left[ u_{j,k}^o - \hat{u}(t_j,x_k;\boldsymbol{q}) \right]^2
\end{equation}
where $\bold{Q}$ is the set of admissible values for the parameters. The simulated and observed values for the population density are denoted as $u_{j,k}^o$ and $\hat{u}(t_j,x_k;\bold{q})$. We would like to emphasize that unlike the traditional inverse problem approach described in this section, the PrMF approach only needs to estimate the weights that describe the probability density function for each pair of $D$ and $\rho$.

\subsection{Data} \label{data}

To evaluate the models' performance, we rely on generated synthetic data, for which we know the true parameters. Specifically, we generate noisy synthetic data by solving the RanDE model (Eqs. \ref{aggregated_discrete_u} and \ref{new_discrete_subpop}). In the previous study, Rutter et al. generated data from two independent double-gaussian distributions for $\boldsymbol{D}$ and $\boldsymbol{\rho}$ \cite{rutter_estimating_2018}. To further generalize the PrMF method, we generate noisy synthetic data from a mixture of gaussian distributions of $\boldsymbol{D}$ and $\boldsymbol{\rho}$. To be consistent with previous studies \cite{hawkins-daarud_quantifying_2019,nardini_learning_2020}, the data is generated over the spatial domain between $x=0$ cm and $x=2$ cm. The medium survival time for GBM patients is about 1.25 years \cite{stupp_radiotherapy_2005}. We extend temporal domain further and we generate the data over $t=[0,1.4]$.  In addition, we choose the number of spatial points, $N_x = 101$, and the number of temporal points, $N_t$ = 51. However, we want to evaluate each model's performance on both fitting and predicting. Hence, we divide the time interval into two sub-intervals: fitting interval ($t=[0,1]$) and predicting interval ($t=[1,1.4]$).

In Table \ref{param_data_table}, we display the means and standard deviations associated with the parameter distributions for $\boldsymbol{D}$ and $\boldsymbol{\rho}$. To incorporate the ``go or grow" hypothesis, we divide the aggregated population into two subpopulations: proliferative population and invasive population. Therefore, we use a mixture of two-gaussian distributions to reflect this ``go or grow" hypothesis. In particular, the proliferative population centers around $(\boldsymbol{D},\boldsymbol{\rho}) = (0.01,10)$, while the invasive population is centered around $(\boldsymbol{D},\boldsymbol{\rho}) = (0.1,1)$. These mean values are chosen to reflect the assumption that proliferative and invasive cells exhibit drastically different behaviors. The invasive cells tend to migrate more compared to the proliferative cells, but the proliferative cells grow faster. Hence, the mean diffusion rate, $\boldsymbol{D}$, for invasive cells is much higher, but the growth rate, $\boldsymbol{\rho}$, is less than the proliferative cells. To test the PrMF method further, we include an additional dataset in which there is an intermediate population with a medium diffusion rate ($\boldsymbol{D} = 0.04$) and a medium growth rate ($\boldsymbol{\rho} = 5$). 

\begin{table}[!htbp]
\caption{Means and standard deviations for the parameters for two and three subpopulations.}\label{param_data_table}%
\begin{tabular}{@{}lccc@{}}
\toprule
\multirow{2}{*}{Distribution} & \multirow{2}{*}{Population} & $\boldsymbol{D}$  & $\boldsymbol{\rho}$ \\
\cmidrule(lr){3-3} \cmidrule(lr){4-4}
& & Mean (Std) & Mean (Std) \\
\midrule
\multirow{2}{*}{Mixture of two-gaussian} & Proliferative/Grow & 0.01 (5e-4) & 10 (1) \\
& Invasive/Go & 0.1 (1e-3) & 1 (1)\\
\midrule
\multirow{3}{*}{Mixture of three-gaussian} & Proliferative/Grow & 0.01 (5e-4) & 10 (1) \\
& Intermediate & 0.04 (5e-4) & 5 (1) \\
& Invasive/Go & 0.1 (1e-3) & 1 (1) \\
\botrule
\end{tabular}
\end{table}

\begin{figure}[!htbp]
\centering
\includegraphics[width=1\linewidth]{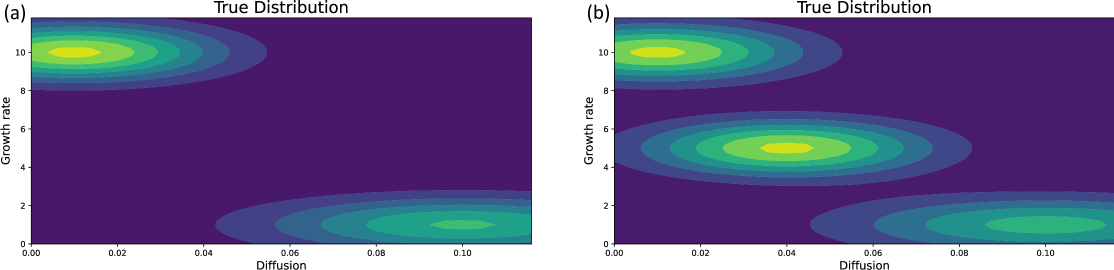}
\caption{ \small True distributions for data generation. (a) shows a mixture of two-gaussian distribution. The top left disc represents the mixture of gaussian for the proliferative population. The bottom right disc represents the mixture of gaussian for the invasive population. Similarly, (b) shows a mixture of three-gaussian distribution with an additionally intermediate population.}
\label{fig:true_distributions}
\end{figure}

For both mixtures of two- and three-gaussian distributions, each subpopulation is generated using \verb|scipy| (a Python package) built-in function, \verb|multivariate_normal|, to generate multivariate normal random variables over the domain $\Omega_D = [0, 0.12]$ and $\Omega_\rho = [0, 12]$ with diagonal covariance matrix, assuming no correlation between $\boldsymbol{D}$ and $\boldsymbol{\rho}$. To compute the overall distribution, we normalize the sum of all subpopulation distributions. In Figure \ref{fig:true_distributions}, we illustrate the true distributions for a mixture of two-gaussian distribution (Figure \ref{fig:true_distributions}a) and a mixture of three-gaussian distribution (Figure \ref{fig:true_distributions}b).

We use central differencing to approximate the spatial derivatives in Eq. \ref{new_discrete_subpop}. Then, we numerically solve the model. After generating the synthetic data, we perturb the aggregated density under the proportional error:
\begin{displaymath}
    u_{j,k}^o = u_{j,k}(1 + \varepsilon_{j,k})
\end{displaymath}
where $\varepsilon_{j,k}$ is normally distributed noise as $\varepsilon \sim \mathcal{N}(\mu,\sigma)$ with $\mu=0$ and $\sigma=0.01$, i.e. $1\%$ proportional noise.
\subsection{Subpopulation clustering}\label{kmean}

After estimating the underlying distribution using PrMF, we use $k$-means clustering, an unsupervised learning method, to discover the number of clusters, i.e. the number of populations from the recovered distribution. We first sample $H$ number of individuals (or realizations) using the estimated joint distribution from the PrMF method. Each individual correspond to $(D^{(i)},\rho^{(i)})$ with $i=1,\dots, H$. The $k$-means clustering method is then applied to the sampled population to group individuals into a pre-defined number of clusters (denoted as $k$ clusters) by computing the distance between each individual and the mean for each cluster. The means are also called the clustered centers for each cluster. Since the number of clusters, $k$, is assumed to be unknown, we perform $k$-means clustering on a range of $k$ values. As the number of clusters increases, the total sum of squares distance between each individual and the means decreases. Therefore, we use the so-called elbow method to help determine the "elbow" or the cutoff for the number of clusters. The elbow is the most optimal number of clusters in the sampled population so adding more clusters does not lower the sum of squares error significantly.

We use the \verb|KMeans| function from the \verb|scikitlearn| package with Lloyd's algorithm using Euclidean distance for clustering. In addition, since $\rho^{(i)}$ are much bigger compared to $D^{(i)}$, we normalize the sampled $(D^{(i)},\rho^{(i)})$ to be within the intervals $[-1,1] \times [-1,1]$ before performing $k$-means clustering. This helps the algorithm avoid biases towards the larger components. Once $k$-means clustering method is performed on multiple values of $k$, we use an elbow plot to determine the number of clusters (subpopulations). We use the \verb|KneeLocator| function from \verb|kneed| package to determine the number of clusters.

\section{Results}\label{results}

In this section, we discuss the results of PrMF in estimating the joint distribution of $\boldsymbol{D}$ and $\boldsymbol{\rho}$ using the delta functions method \cite{rutter_estimating_2018,banks_comparison_2007}. As mentioned in Section \ref{forwardsolve_PMF}, we first needed to forward solve for the phenotype cell densities, $c(x,t, D_i,\rho_i)$ over $M = M_D \times M_\rho$ finely meshed nodes for $\boldsymbol{D}$ and $\boldsymbol{\rho}$. As the number of nodes increases, the PrMF method could become prone to overfitting as it might attempt to fit the noise. Therefore, it is necessary to find the most optimal model that generates the best fit. In this study, we tested the PrMF method over different combinations of $M_D \times M_\rho$. We then used the Akaike Information Criteria (AIC) \cite{akaike_new_1974,banks_aic_2017} to penalize the models with the higher number of nodes and find the most optimal model that generates the best fit. 

In this study, we chose $\max_{M_D} = \max_{M_\rho} = 20$. For a more computationally efficient approach, we only forwarded solve over a $\max_{M_D} \times  \max_{M_\rho}$ mesh. When testing the PrMF method over different combinations of $M_D \times M_\rho$ (with $M_D = 5, 10, 20$ and $M_\rho =5,10, 20$), we generated lower-resolution meshes with $M_D \times M_\rho$ from $\max_{M_D} \times  \max_{M_\rho}$ and perform PrMF on the lower-resolution meshes. Afterward, we performed $k$-means clustering on the recovered mesh with the lowest AIC score. In this section, we only show the results for generated data set with a mixture of two-gaussian distribution. The results for a mixture of three-gaussian distribution are in the Supplemental Materials \ref{SM_3P}.
\subsection{Inverse problem result using Prokhorov metric framework}

\subsubsection{Fitting and forecasting results}
\label{fitting_forecasting_section}
\begin{figure}[!htbp]
\centering
\includegraphics[width=1\linewidth]{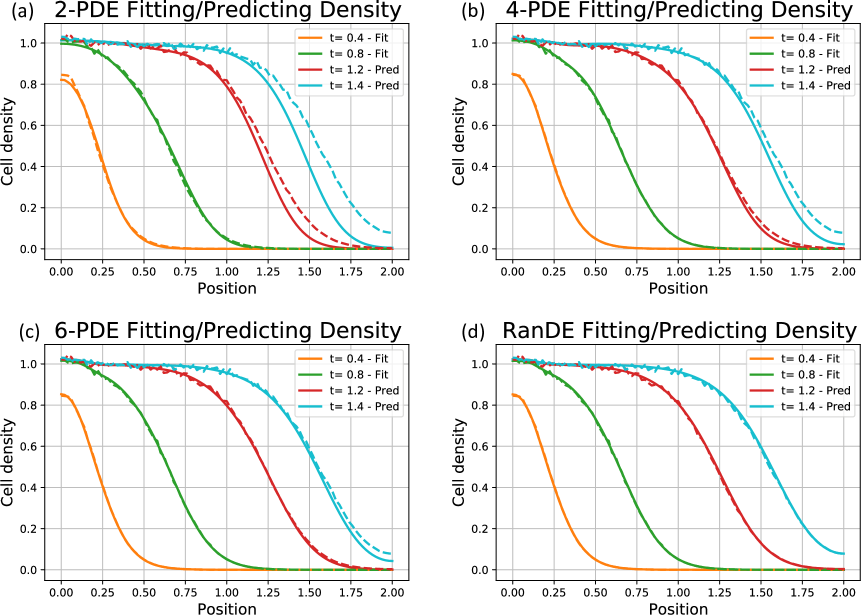}
\caption{ \small Aggregated cell density comparison between: (a) 2-PDE model, (b) 4-PDE model, (c) 6-PDE model, and (d) RanDE model. In each figure, we plot the generated data (dashed curves) and model simulation (solid curves) for 4 different time points. For fitting interval, we plot the cell density at $t=0.4$ and $t=0.8$. For the prediction interval, we plot the cell density at $t=1.2$ and $t=1.4$.}
\label{fig:fitting_forecasting_results}
\end{figure}

In this section, we compare RanDE performance against the traditional inverse problem method on models with 2-PDE, 4-PDE, and 6-PDE. The 2-PDE, 4-PDE, and 6-PDE models assume a known subpopulation of 2, 4, and 6, respectively. Fig. \ref{fig:fitting_forecasting_results} illustrates the simulated results on data with the mixture of two-gaussian, using the estimated parameters from the traditional inverse problem approach on 2-PDE, 4-PDE, and 6-PDE (Figs. \ref{fig:fitting_forecasting_results}a-\ref{fig:fitting_forecasting_results}c) and the simulated aggregated population from the estimated RanDE model using the PrMF approach (Fig. \ref{fig:fitting_forecasting_results}d). In Fig. \ref{fig:fitting_forecasting_results}, we plot the cell densities for two representative temporal time points within the fitting interval, $t = 0.4$ and $0.8$, and another two within the prediction time interval, $t= 1.2$ and $1.4$. The results in Fig. \ref{fig:fitting_forecasting_results} show that all models were able to fit the generated data. However, in terms of forecasting, the RanDE model is more capable of describing the unseen densities at $t=1.2$ and $t=1.4$. In fact, as we increase the number of PDEs in the non-RanDE models from 2 to 6, the prediction results tend to converge to the RanDE prediction result. Similar results for the mixture of three-gaussian can be found in Fig. \ref{fig:fitting_forecasting_results_3P}.

\begin{figure}[!htbp]
\centering
\includegraphics[width=1\linewidth]{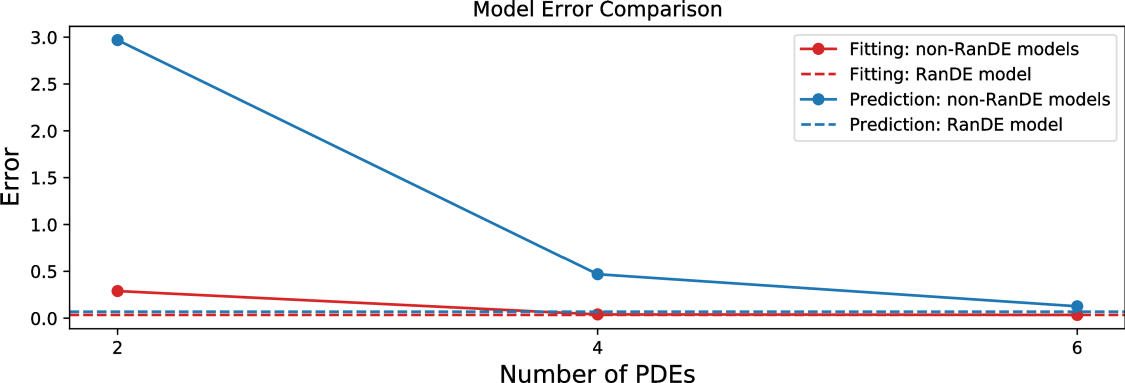}
\caption{ \small Comparison for the fitting and prediction errors between models. The solid red and blue curves show the $\textit{SSE}$ for non-RanDE models with the number of on the x-axis within fitting and prediction intervals, respectively. The red and blue horizontal dashed lines are the $\textit{SSE}$ for the RanDE model within fitting and prediction intervals, respectively. }
\label{fig:fitting_forecasting_error}
\end{figure}

In Fig. \ref{fig:fitting_forecasting_error}, we plot the $\textit{SSE}$ comparison between models. We found the 2-PDE model has a higher fitting error compared to other models, however, the differences are not significant. Additionally, while the RanDE model has a similar fitting error compared to the 4-PDE and 6-PDE models, the RanDE model is better in terms of forecasting. We obtained similar results for the mixture of three-gaussian (See Fig. \ref{fig:fitting_forecasting_error_3P}).

\subsubsection{Profile of traveling wave speed}
\label{wave_speed_profile_section}
One approach to measuring the tumor spread is to use the traveling wave speed. Due to the complexity of the models, in this work, we computed the traveling wave speed numerically \cite{stepien_data-motivated_2015,stepien_traveling_2018}.  To compute the traveling wave speed, we measure the position $x$, where the aggregated cell density equals a specified density value, $u^*$, for different time values $t$. To ensure that the wave profile is established for the fitting interval, we ignored the first few time points and only computed the wave speed within the time interval $t=[0.6,1]$. Since the wave shape could potentially change due to the heterogeneity in the cell population, we compute the wave speed for different cell density values ranging from $0.25$ to $0.65$.

\begin{figure}[!htbp]
\centering
\includegraphics[width=1\linewidth]{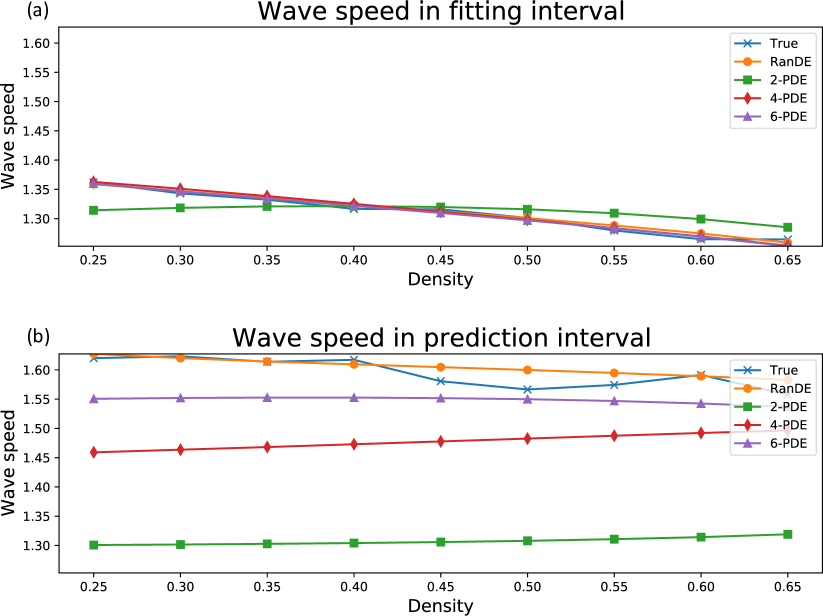}\caption{ \small Wave speed profile comparison between models within: (a) fitting interval and (b) prediction interval.}
\label{fig:wave_speed}
\end{figure}

In this section, we compare the computed profiles for the traveling wave speed at different cell densities between the models on data with mixture of two-gaussian (for mixture of three-gaussian, see \ref{SM_wave_speed_3P}). In Fig. \ref{fig:wave_speed}, we plot the computed traveling wave speed profiles within the fitting interval (Fig. \ref{fig:wave_speed}a) and within the forecasting interval (Fig. \ref{fig:wave_speed}b) from the data and the models (See Fig. \ref{fig:wave_speed_3P} for mixture of three-gaussian). In Fig. \ref{fig:wave_speed}a, we found that the computed profiles for wave speed from all models, except the profile from the 2-PDE model, resemble the computed profile from the generated data during the fitting interval. However, for the forecasting interval, only the profile from the RanDE model is able to approximate the true profile. This reinforced the observation in Section \ref{fitting_forecasting_section}, in which, we found that while the 4-PDE and 6-PDE models are able to fit the data, their forecasting performances are less accurate than the RanDE model.

\subsubsection{Recovered distribution and cluster centers}
\label{recovered_joint_distribution_section}
\begin{figure}[!htbp]
\centering
\includegraphics[width=1\linewidth]{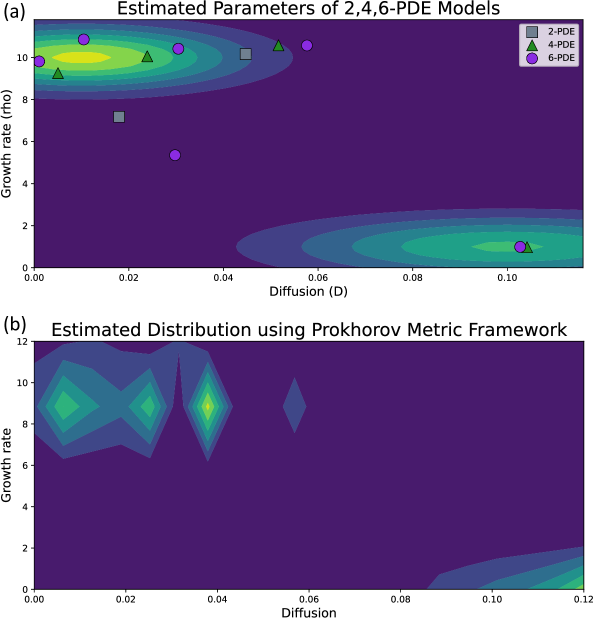}
\caption{ \small (a) Point-wise estimated parameters of the 2,4,6-PDE models on the true distribution with 30 $D$-nodes and 60 $\rho$-nodes. (b) Estimated distribution using Prokhorov metric framework with 20 $D$-nodes and 5 $\rho$-nodes.}
\label{fig:param_mesh}
\end{figure}

After showing that the RanDE model is capable of fitting and forecasting in Section \ref{fitting_forecasting_section}, we then investigate the PrMF method's ability to recover the underlying distribution from the noisy synthetic data generated from mixture of two-gaussian (See Appendix \ref{recovered_dist_3P} for mixture of three-gaussian). In Fig. \ref{fig:param_mesh}, we compare the true distribution (Fig. \ref{fig:param_mesh}a) and the estimated distribution using PrMF (Fig. \ref{fig:param_mesh}b). While the PrMF estimated distribution is coarser compared to the true distribution, the PrMF method is able to identify two distinct clusters of densities, $P(\boldsymbol{D},\boldsymbol{\rho})$, as shown in the true distribution. Furthermore, we plot the point-wise estimated values of $(D,\rho)$ for the non-RanDE models in Fig. \ref{fig:param_mesh}a. Interestingly, these non-RanDE models focus more on estimating the proliferative subpopulation, especially the 2-PDE model.

\begin{figure}[!htbp]
\centering
\includegraphics[width=1\linewidth]{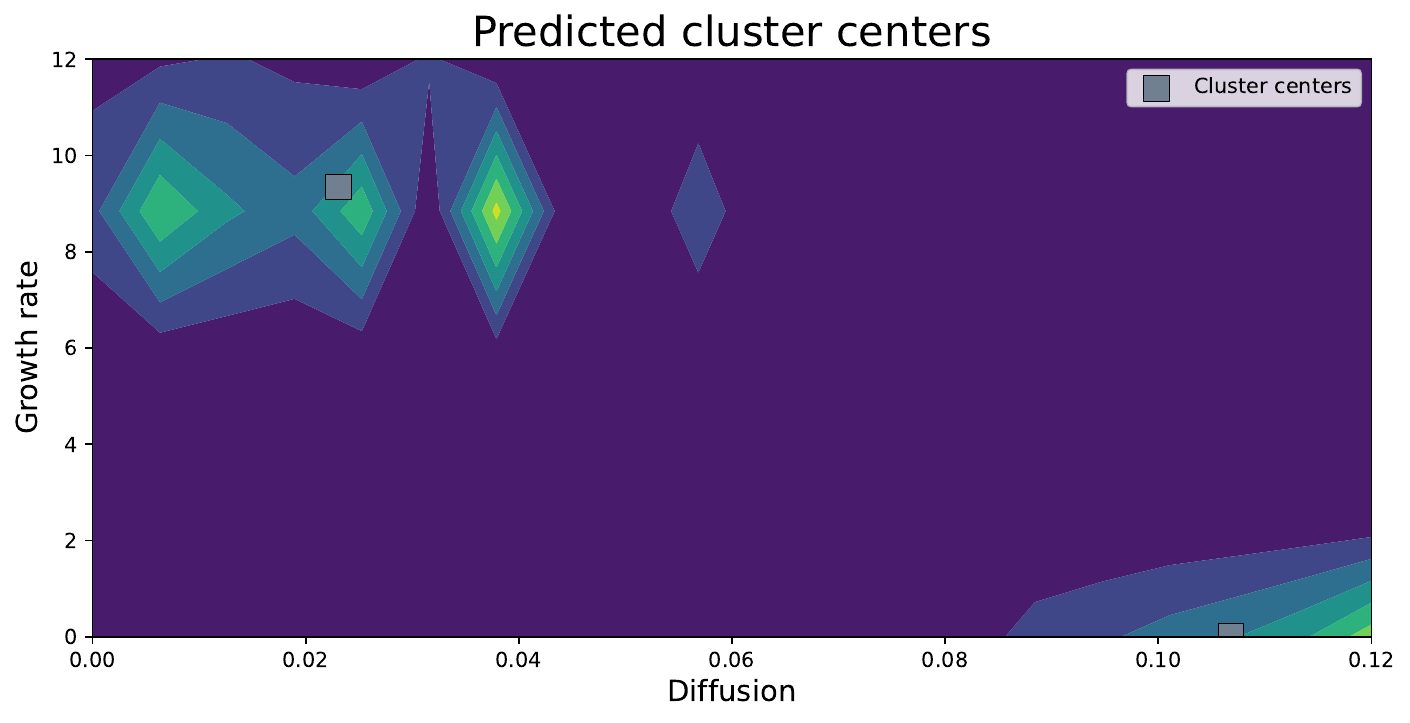}
\caption{ \small Plotting the predicted cluster centers using $k$-means clustering. } \label{fig:cluster_centers}
\end{figure}

We showed that PrMF is able to recover the underlying distribution of the generated data. Our next step is to identify the number of subpopulations (or clusters) and the cluster centers from the recovered distribution. We first sampled $H = 10,000$ individual cells from the recovered distribution. Then, we used $k$-means clustering to group sampled individual cells into $k$ clusters, with $k$ varying from 1 to 10. Using the elbow plot (see Fig. \ref{fig:2P_elbow}), we found that $k=2$ is the most optimal number of clusters. This means there are 2 subpopulations that can be generated from the recovered distribution. This agrees with the fact that we generated the aggregated data from 2 subpopulations. In Fig. \ref{fig:cluster_centers}, we plot the cluster centers that were identified by $k-$means clustering. We note that, unlike the results from the 2-PDE model, our cluster centers were able to accurately identify both distinct subpopulations: the ``go" subpopulation with high $D$, low $\rho$ and the ``grow" subpopulation with low $D$, high $\rho$. When we used the 2-PDE model which has 2 subpopulations, we recovered two ``grow" populations that had high growth rates and low diffusion rates.

\subsubsection{Computational cost}
\label{computational_cost_section}

\begin{figure}[!htbp]
\centering
\includegraphics[width=1\linewidth]{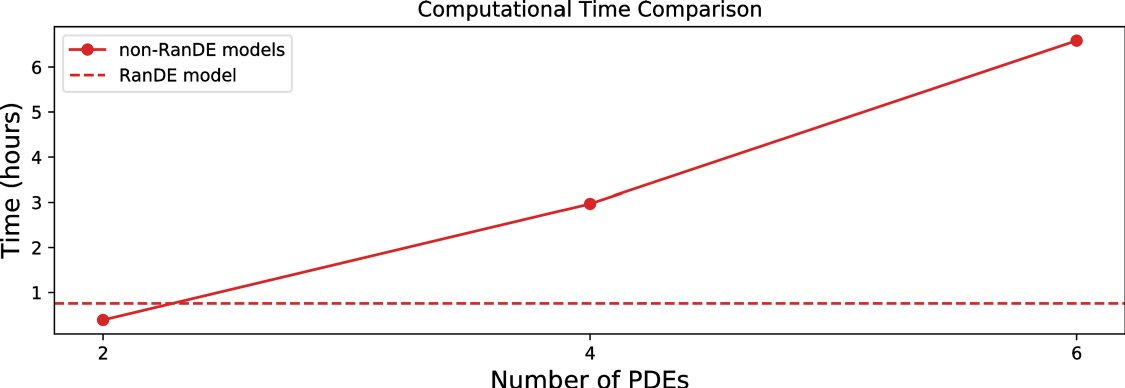}
\caption{ \small Comparison of the computational time between models. The solid curve represents the computational time required to perform traditional inverse problems on non-RanDE models. The horizontal dashed line represents the computational time required to perform the inverse problem using Prokhorov metric framework on the RanDE model. }
\label{fig:computational_cost}
\end{figure}

In this section, we provide details on the computational cost for the inverse problem on each model. The optimization problems for all 4 models are non-linear optimization. Therefore, we performed optimization for each model at 20 different randomly chosen starting sets of parameters and only record the optimized parameter sets with the lowest $\textit{SSE}$. However, the inverse problem approach using PrMF required additional steps prior to the optimization step. In particular, we needed to forward solve for the solution $c(x,t; D_i,\rho_i)$ of Eq. \ref{new_discrete_subpop_data} for different sets of $(D_i,\rho_i)$. In this work, we chose $\max_{M_D} = \max_{M_\rho} = 20$. Therefore, 400 forward solved solutions, $c(x,t;D_i,\rho_i)$, were computed. In addition, for each combination of $M_D$ and $M_\rho$, we performed optimization at 20 different randomly chosen starting sets of probability measures. In Fig. \ref{fig:computational_cost}, we plot the computational time between models in hours. For the traditional inverse problem approach on non-RanDE models, as we increased the number of PDEs, the computational time for the inverse problem increased. Additionally, performing inverse problems using PrMF on the RanDE model required more time compared the traditional approach on the 2-PDE model. However, it is much more time-efficient to perform inverse problem using PrMF on RanDE model than the traditional approach on the 4-PDE and 6-PDE models.

\section{Discussion and conclusions}\label{discussion_section}
Building upon previous work \cite{rutter_estimating_2018}, we have introduced a new RanDE model that is an extension of the ``go or grow" hypothesis on GBM cells. In addition, we investigated the ability of the PrMF method on recovering the underlying distribution from the aggregated cell density data that we generated from the RanDE model. However, our work has addressed two limitations of the previous work. First, unlike the previous work \cite{rutter_estimating_2018}, the new random differential equation model incorporates the intratumor competition between phenotypes. This  allows the proposed model to be more biologically realistic compared to the previous model. Second, rather than using the assumption of independence between parameter distributions, we assumed a joint probability density distribution. This assumption allows the model to generalize the ``go or grow" hypothesis in which one subpopulation exhibits high proliferation and slow diffusion while another type of subpopulation exhibits low proliferation and high diffusion.

Our results for mixture of two-gaussian show that the RanDE model outperforms the 2-PDE model with competition in both fitting and forecasting (see Figs. \ref{fig:fitting_forecasting_results} and \ref{fig:fitting_forecasting_error}). In addition, while the RanDE has similar performance in terms of fitting compared to the 4-PDE and 6-PDE models, it is capable of forecasting future cell density. Furthermore, using the RanDE model simulation, we were able to compute much more accurate wave speed profiles in the forecasting interval (see Fig. \ref{fig:wave_speed}). In fact, as we increase the number of coupled PDEs in the non-RanDE models from 2 to 6, the prediction results tend to converge to the RanDE prediction result. This can be indirectly observed in Fig. \ref{fig:fitting_forecasting_results}a-\ref{fig:fitting_forecasting_results}d, in which, the non-RanDE models cell density curves within the prediction interval ($t=1.2$ and $t=1.4$) start to get closer to the true cell density curves. On the other hand, the predicted cell densities from the RanDE model are indistinguishable from the true cell densities. Fig. \ref{fig:fitting_forecasting_error} further confirms this observation. As we increase the number of coupled PDEs from 2 to 6, the prediction error for non-RanDE models decreases to approach the RanDE model prediction error. These results lead to the computed wave speed using the coupled non-RanDE models getting closer to the computed wave speed using the RanDE model as we increase the number of PDEs from 2 to 6  (see Fig. \ref{fig:wave_speed}b). We obtained results for data generated using a mixture of three-gaussian (see Figs. \ref{fig:fitting_forecasting_results_3P}, \ref{fig:fitting_forecasting_error_3P}, \ref{fig:wave_speed_3P}).

In Fig. \ref{fig:param_mesh}b, we showed that PrMF is able to recover the underlying mixture of two-gaussian distribution. We also demonstrated that with a simple unsupervised machine learning method ($k$-means clustering), we were able to predict the number of subpopulations within the aggregated cell density data. In addition, we found that performing inverse problem using PrMF on the RanDE model is more time efficient compared to the traditional inverse problem on 4-PDE and 6-PDE models. Similar results for mixture of three-gaussian can be found in Fig. \ref{fig:param_mesh_3P}b. However, $k$-means clustering only predict two clusters of individuals rather than three clusters. This could be due to the fact that recovered distribution is much coarser compared to the true distribution.

While the proposed RanDE model is much more complicated in terms of dimensionality compared to other PDE models, it is still lacking other details about the GBM growth and proliferation. One such example is to include an advection term \cite{stein_mathematical_2007}. Therefore, this work should be tested further for three or more mixtures of gaussian distributions. 

One current disadvantage of the PrMF method is that it requires a good approximation on the boundaries of the parameters in order to perform the forward-solving step. For future work, optimal design methods such as \emph{SE}-optimal design \cite{banks_comparison_2011,banks_experimental_2013,adoteye_optimal_2015} should be considered to help generate an adaptive mesh to perform forward solving. \emph{SE}-optimal design has been applied to find the optimal observation time to maximize the information gained from the experimental data \cite{banks_comparison_2011,adoteye_optimal_2015}.

In future work, the accuracy of the PrMF applied to the RanDE model can be evaluated in combination with data denoising methods. For example, in \cite{lagergren_learning_2020} the authors showed that data arising from a reaction diffusion model, similar to the data we investigated in this work, could be accurately approximated using a neural network. Importantly, the neural network approximated the time and space derivatives up to second order more accurately than finite difference or spline based methods. In addition, we postulate that the PrMF method can be combined with model-free statistical error model inference methods, such as the method developed in \cite{banks_use_2016}. We note that the work presented here used the PrMF with an ordinary least squares cost function, even though the noisy data were generated according to a proportional error model that would be more compatible with a generalized least squares cost function \cite{nardini_learning_2020}. Thus, we did not assume to know whether the noise was generated by a proportional or constant error observation process. Methods such as those described in \cite{banks_use_2016} can be used to infer the most appropriate statistical error model, i.e., proportional or constant error, and to approximate the variance of the residuals all without requiring a mathematical model. In practice, such data denoising and statistical model inference methods could both be applied prior to using the PrMF to estimate parameter distributions, and the PrMF could be adapted to use a generalized least squares framework \cite{banks_prohorov_2018}.

Finally, the PrMF method can be applied to any modeling problem to estimate both point-wise parameters or the distribution of parameters. It has been previously used on other biological problems \cite{banks_estimation_1999, banks_comparison_2007,banks_estimating_2018,schacht_estimation_2019,banks_parameter_2020}. Therefore, we plan to publish the codes for anyone who is interested in using PrMF on their modeling problem. 
\backmatter
\section*{Declarations}
\bmhead{Funding}
Kyle Nguyen was supported by the National Science Foundation Graduate Research Fellowship under Grant No. DGE-2137100.

\bmhead{Acknowledgments}
We would like to thank Celia Schacht for her helpful comments.

\bmhead{Data availability statement}
The code and data are publicly available at: \url{https://github.com/kcnguyen3191/rande_prmf}

\begin{appendices}
\section{Mixture of two-gaussian distribution}
\subsection{Elbow plot} \label{2P_eblow}
\begin{figure}[!htbp]
\centering
\includegraphics[width=1\linewidth]{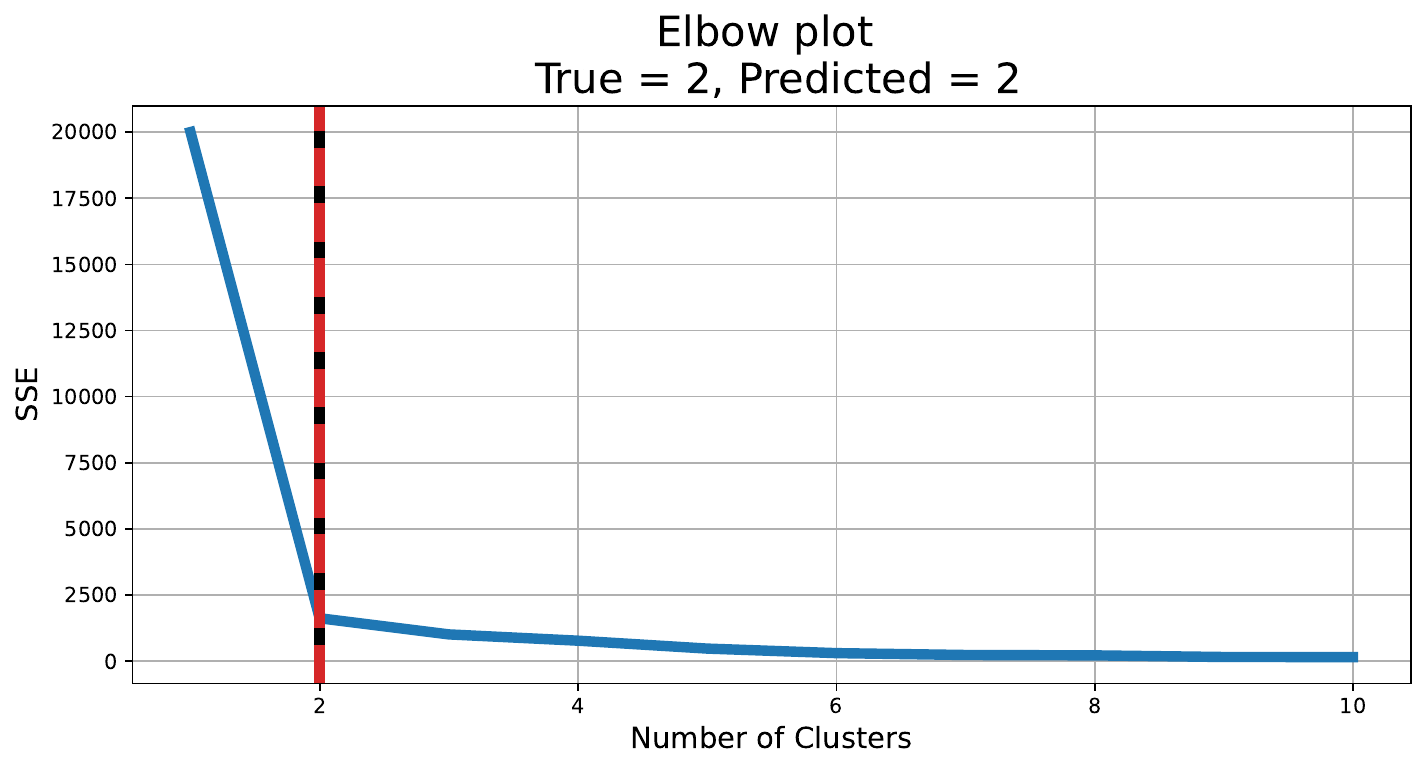}
\caption{ \small Elbow plot showing the sum of squares error for different values of $k$.} \label{fig:2P_elbow}
\end{figure}

Fig. \ref{fig:2P_elbow} displays the elbow plot showing the sum of squares error for different values of $k$. The true number of subpopulations (2) is shown in the black solid line. \verb|KneeLocation| is then used to determine the number of cluster centers, which is shown in the red dashed line.

\section{Mixture of three-gaussian distribution} \label{SM_3P}
\subsection{Fitting and forecasting results}
\begin{figure}[!htbp]
\centering
\includegraphics[width=1\linewidth]{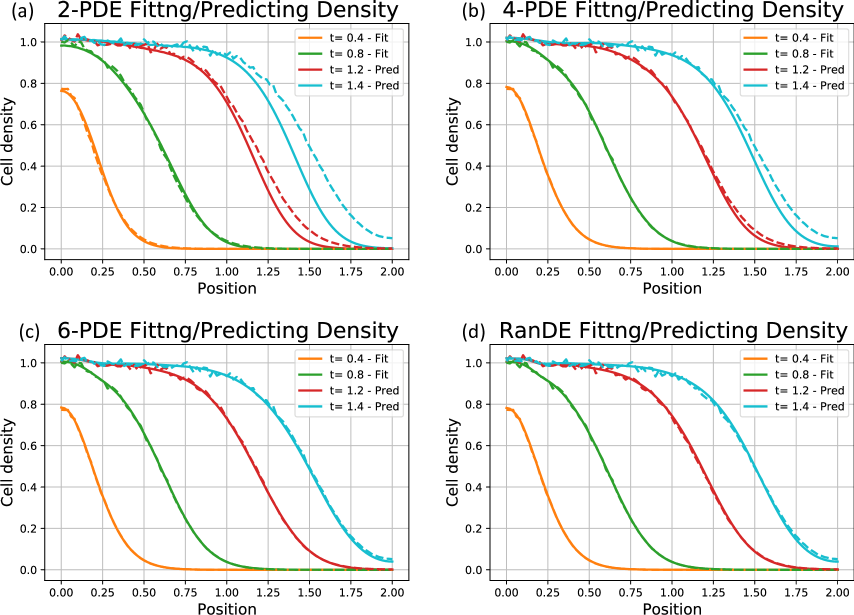}
\caption{ \small Aggregated cell density comparison between: (a) 2-PDE model, (b) 4-PDE model, (c) 6-PDE model, and (d) RanDE model. In each figure, we plot the generated data (dashed curves) and model simulation (solid curves) for 4 different time points. For fitting interval, we plot the cell density at $t=0.4$ and $t=0.8$. For the forecasting interval, we plot the cell density at $t=1.2$ and $t=1.4$.}
\label{fig:fitting_forecasting_results_3P}
\end{figure}

\begin{figure}[!htbp]
\centering
\includegraphics[width=1\linewidth]{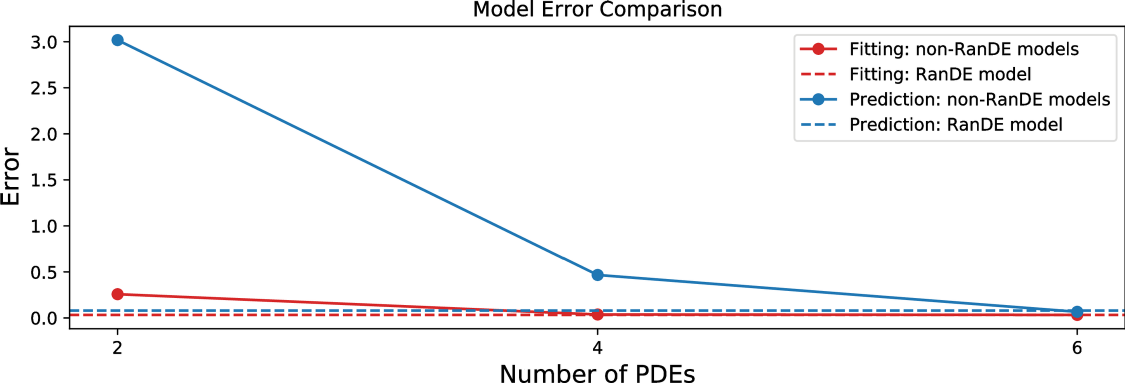}
\caption{ \small Comparison for the fitting and prediction errors between models. The solid red and blue curves show the $\textit{SSE}$ for non-RanDE models with the number of on the x-axis within fitting and prediction intervals, respectively. The red and blue horizontal dashed lines are the $\textit{SSE}$ for the RanDE model within fitting and prediction intervals, respectively. }
\label{fig:fitting_forecasting_error_3P}
\end{figure}

In Fig. \ref{fig:fitting_forecasting_results_3P}, we plot the simulated results using the estimated parameters from the traditional inverse problem approach on 2-PDE, 4-PDE, and 6-PDE (Figs. \ref{fig:fitting_forecasting_results_3P}a-\ref{fig:fitting_forecasting_results_3P}c) and the simulated aggregated population from the estimated RanDE model using PrMF approach (Fig. \ref{fig:fitting_forecasting_results_3P}d). In Fig. \ref{fig:fitting_forecasting_error_3P}, we plot the $\textit{SSE}$ comparison between models.

\subsection{Profile of traveling wave speed} \label{SM_wave_speed_3P}
\begin{figure}[!htbp]
\centering
\includegraphics[width=1\linewidth]{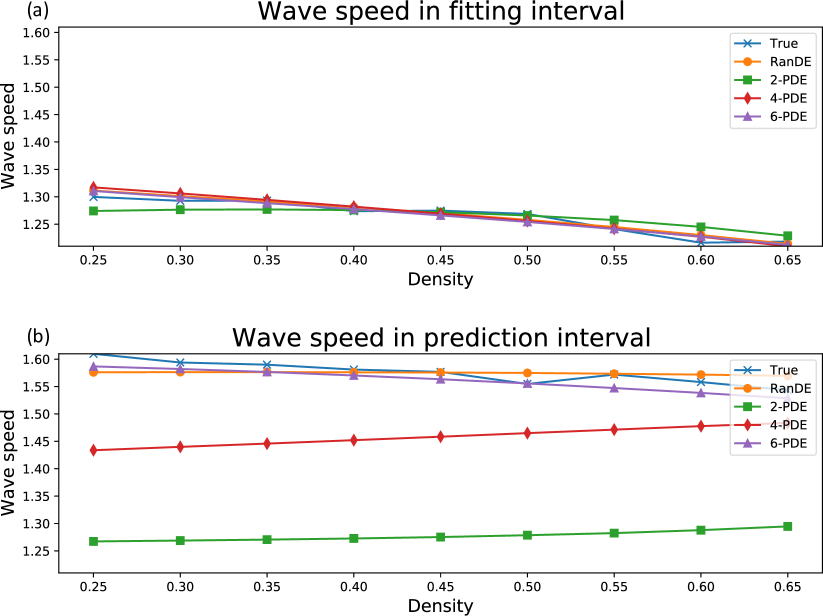}
\caption{ \small Wave speed profile comparison between models within: (a) fitting interval and (b) prediction interval.}
\label{fig:wave_speed_3P}
\end{figure}

In Fig. \ref{fig:wave_speed_3P}, we plot the estimated traveling wave speed within the fitting interval (Fig. \ref{fig:wave_speed_3P}a) and within the forecasting interval (Fig. \ref{fig:wave_speed_3P}b).

\subsection{Recovered distribution and cluster centers} \label{recovered_dist_3P}
\begin{figure}[!htbp]
\centering
\includegraphics[width=1\linewidth]{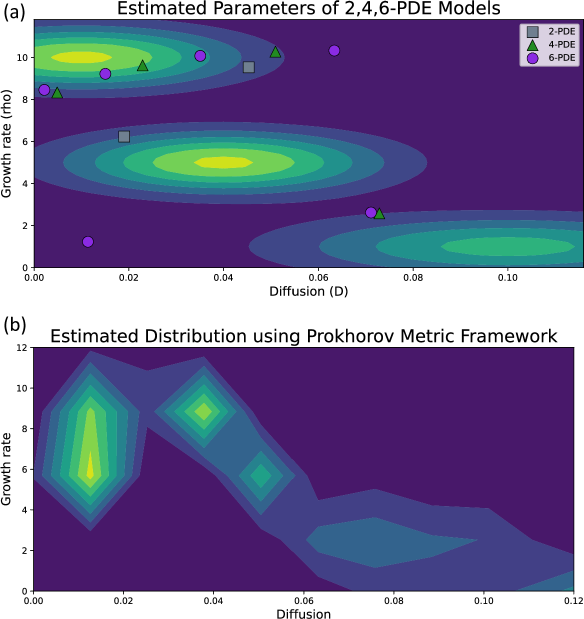}
\caption{ \small Point-wise estimated parameters of the 2,4,6-PDE models on the true distribution with 30 $D$-nodes and 60 $\rho$-nodes. (b) Estimated distribution using Prokhorov metric framework with 10 $D$-nodes and 5 $\rho$-nodes.}
\label{fig:param_mesh_3P}
\end{figure}

\begin{figure}[!htbp]
\centering
\includegraphics[width=1\linewidth]{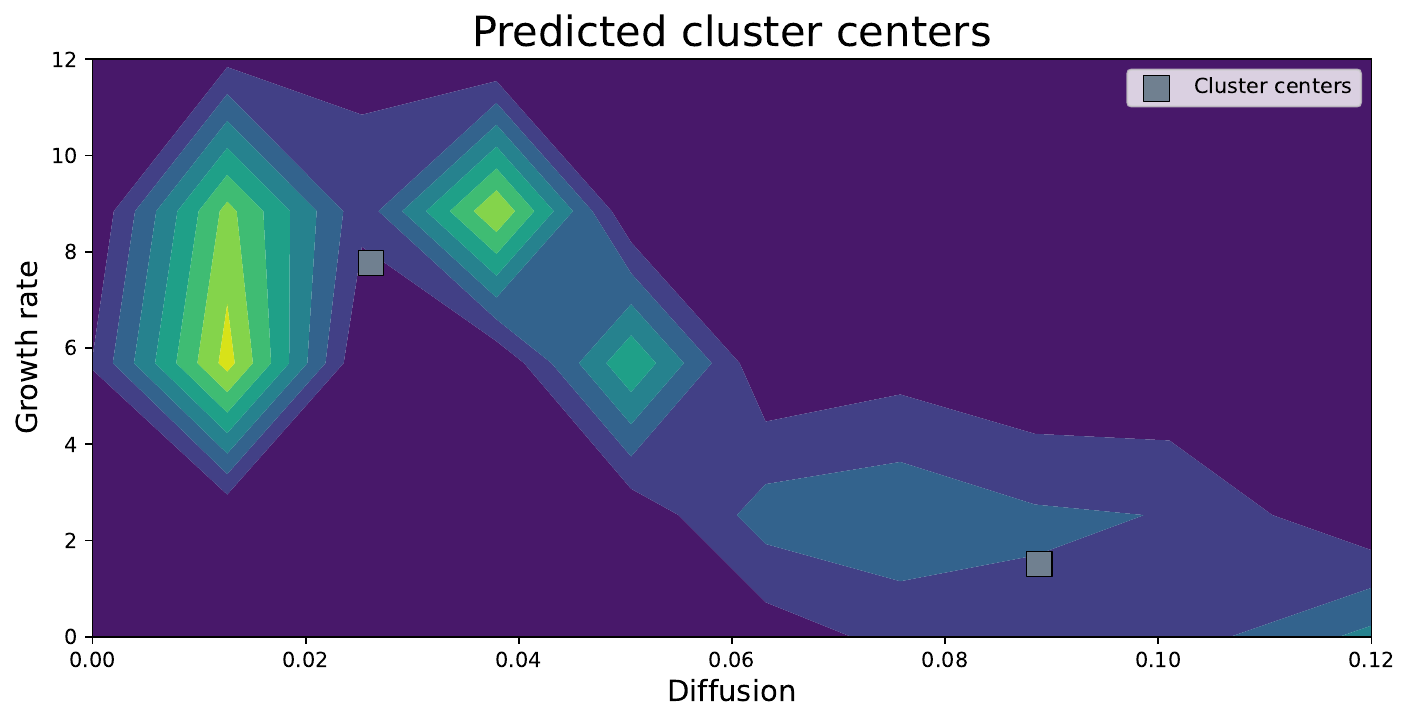}
\caption{ \small Plotting the predicted cluster centers using $k$-means clustering. } \label{fig:cluster_centers_3P}
\end{figure}

In Fig. \ref{fig:param_mesh_3P}, we compare the true distribution (Fig. \ref{fig:param_mesh_3P}a) and the estimated distribution using PrMF (Fig. \ref{fig:param_mesh_3P}b). In Fig. \ref{fig:cluster_centers_3P}, we plot the cluster centers that were identified by $k$-means clustering.

\subsection{Elbow plot} \label{3P_eblow}
\begin{figure}[!htbp]
\centering
\includegraphics[width=1\linewidth]{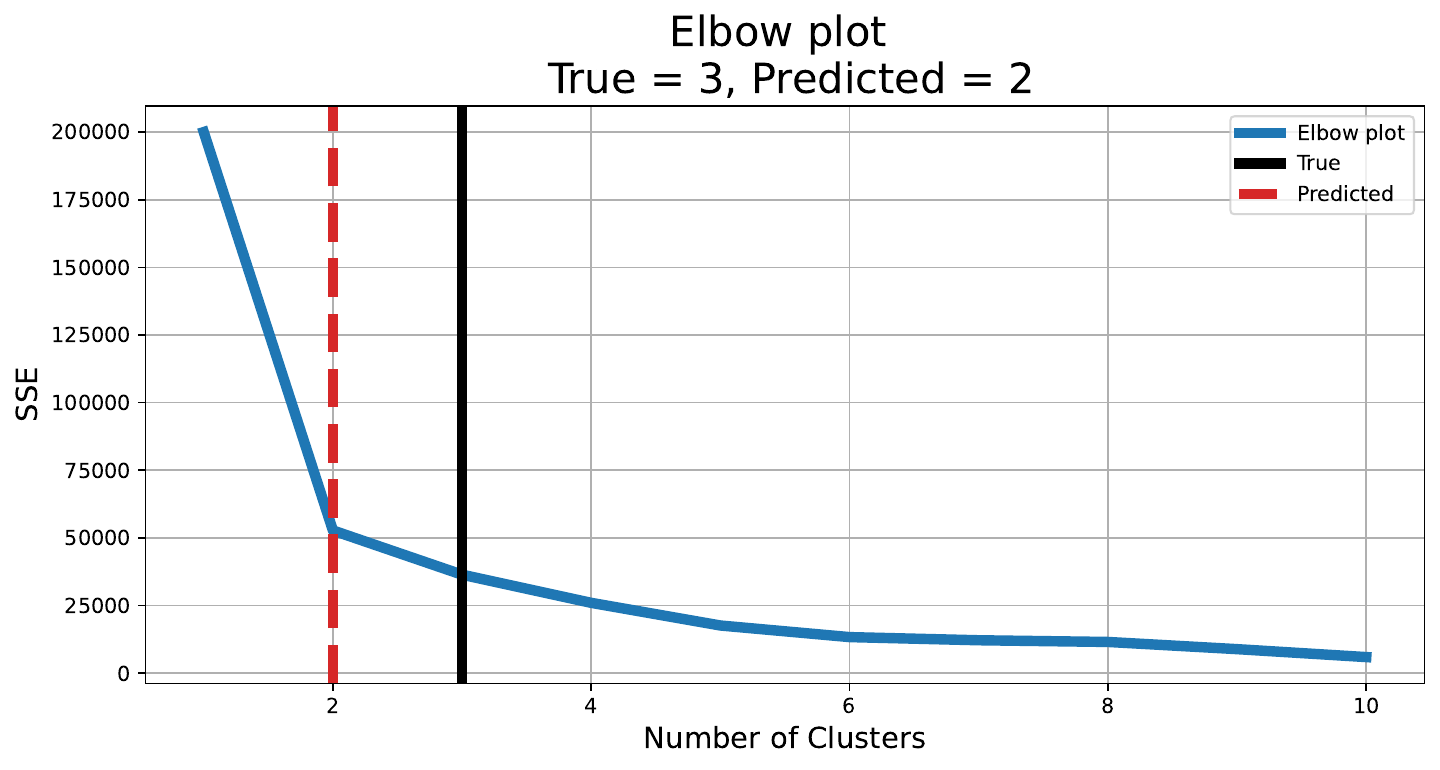}
\caption{ \small Elbow plot showing the sum of squares error for different values of $k$.}
\label{fig:3P_elbow}
\end{figure}

In Fig. \ref{fig:3P_elbow}, we plot the elbow curve sum of squares error curve against the number of clusters, $k$. We find that the predicted number of clusters in the elbow plot (red dashed line) is two.

\end{appendices}


\bibliography{sn-bibliography}

\end{document}